\documentclass{elsart}
\usepackage{mathrsfs}
\usepackage{epsfig,latexsym,amsfonts,amsmath,amssymb}

\begin{document}
\begin{frontmatter}

\newtheorem{theorem}[subsection]{Theorem}
\newtheorem{lemma}[subsection]{Lemma}
\newtheorem{conjecture}[subsection]{Conjecture}
\newtheorem{proposition}[subsection]{Proposition}
\newtheorem{definition}[subsection]{Definition}
\newtheorem{corollary}[subsection]{Corollary}
\newtheorem{example}[subsection]{Example}
\newtheorem{remark}[subsection]{Remark}
\newtheorem{de}[subsection]{Definition}
\newtheorem{problem}[subsection]{Problem}

\renewcommand{\theequation}{\arabic{section}. \arabic{equation}}
\renewcommand{\thesection}{\arabic{section}}
\renewcommand{\thethm}{\arabic{section}.\arabic{thm}}
\renewcommand{\thefootnote}{\fnsymbol{footnote}}

\newcommand{\df}{\stackrel{\mbox{\rm def }}{=}}
\newcommand{\di}{\displaystyle}
\newcommand{\bl}[2]{{\left\langle #1 \:\:  \vrule \:\:  #2 \right\rangle}}
\newcommand{\vanish}[1]{}

\title{Independent Sets in Direct Products of Vertex-transitive Graphs}

\author{Huajun  Zhang}\ead{huajunzhang@zjnu.cn}

\address{ Department of Mathematics,
 Zhejiang Normal University, Jinhua 321004, P.R. China}

\address{ Department of Mathematics, Shanghai Normal University,
Shanghai 200234,  P.R. China}
\baselineskip 20pt

\date{}
\maketitle
 \vspace{4mm}

\begin{abstract}
The direct product $G\times H$ of graphs $G$ and $H$ is defined by:
\[V(G\times H)=V(G)\times V(H)\]
and
\[E(G\times H)=\left\{[(u_1,v_1),(u_2,v_2)]: (u_1,u_2)\in E(G) \mbox{\ and\ } (v_1,v_2)\in E(H)\right\}.\]
 In this paper, we will prove that the equality
$$\alpha(G\times H)=\max\{\alpha(G)|H|, \alpha(H)|G|\}$$
holds for all vertex-transitive graphs $G$ and $H$, which provides an affirmative answer to a problem posed by Tardif (Discrete Math. 185 (1998) 193-200). Furthermore, the structure of all maximum independent sets of $G\times H$ are determined.

\begin{keyword} direct product; primitivity; independence number; vertex-transitive
 \\[7pt]
\end{keyword}
\end{abstract}
\end{frontmatter}

\newcommand{\lr}[1]{\langle #1\rangle}
\newcommand{\qchoose}[2]{{ #1   \atopwithdelims[]  #2 }}

\parindent 17pt
\baselineskip 17pt

\section{Introduction}
Let $G$ and $H$ be two graphs. The direct product $G\times H$ of $G$ and $H$ is defined by£º
\begin{eqnarray*}
V(G\times H)=V(G)\times V(H)\end{eqnarray*}
and \begin{eqnarray*}
E(G\times H)=\left\{[(u_1,v_1),(u_2,v_2)]: (u_1,u_2)\in E(G) \mbox{\ and\ } (v_1,v_2)\in E(H)\right\}.
\end{eqnarray*}
It is easy to see this  product is commutative and associative, and  the  product
of more than two graphs is well-defined. For a graph $G$, the products $G^n=G\times G\times \cdots\times G$ is called
the $n$-th powers of $G$.

An interesting problem is the independence number of $G\times H$. It is clear that if $I$ is an independent
 set of $G$ or $H$, then the preimage of $I$ under projections is an independent set of $G\times H$, and so
$\alpha(G\times H)\geq \max\{\alpha(G)|H|, \alpha(H)|G|\}.$ It is natural to ask whether the equality holds or not.
In general, the equality does not hold for non-vertex-transitive graphs (see
\cite{Kla}). So Tardif \cite{Tardif} posed the following problem.
\begin{problem}\label{Tard}(Tardif \cite{Tardif})
Does the equality \[\alpha(G\times H)=\max\{\alpha(G)|H|,\alpha(H)|G|\}\] hold for all vertex-transitive
graphs $G$ and $H$?
\end{problem}

Furthermore, it immediately raises another interesting problem:
\begin{problem}\label{MIS}
When $\alpha(G\times H)=\max\{\alpha(G)|H|,\alpha(H)|G|\}$, is every
maximum independent set of $G\times H$  the preimage of an
independent set of one factor under projections?
\end{problem}
If the answer is yes, we then say the direct product $G\times H$ is
MIS-normal (maximum-independent-set-normal). Furthermore, the direct
products $G_1\times G_2\times\cdots\times G_n$ is said to be
MIS-normal if every maximum independent set of it is the preimage of
an independent set of one factor under projections.

About these two problems, there are some progresses have been made for some very special vertex-transitive graphs.

 Let $n,r$ and $t$ be three integers with $n\geq r\geq t\geq 1$. The graph $K(t,r,n)$ is defined by: whose vertices set is the set of all $r$-element subsets of $[n]=\{1,2,\ldots,n\}$, and $A$ and $B$ of which  are adjacent if
and only if $|A\cap B|<t$. If $n\geq 2r$, then $K(1,r,n)$ is the
well-known Kneser graph. The classical Erd\H{o}s-Ko-Rado Theorem
\cite{EKR} states that $\alpha(K(1,r,n))=\binom{n-1}{r-1}$ (where
$n\geq 2r$), and Frankl \cite{frankl3} first investigated the
independence number of the direct products of Kneser graphs.
Subsequently, Ahlswede, Aydinian and Khachatrian investigated the
general case \cite{AAK}.
\begin{theorem}
Let $n_i\geq r_i\geq t_i$ for $i=1,2,\ldots, k$. \\
(i) (Frankl \cite{frankl3}) if $t_1=\cdots=t_k=1$ and $\frac{r_i}{n_i}\geq \frac{1}{2}$ for  $i=1,2,\ldots,k$, then $$\alpha\left(\prod_{1\leq i\leq
k}K(1,r_i,n_i)\right)=\max\left\{\frac{r_1}{n_1},\frac{r_2}{n_2},\ldots,\frac{r_k}{n_k}\right\}\prod_{1\leq
i\leq k}|K(1,r_i,n_i)|.$$\\
(ii) (Ahlswede, Aydinian and Khachatrian \cite{AAK}) $$\alpha\left(\prod_{1\leq i\leq
k}K(t_i,r_i,n_i)\right)=\max\left\{\frac{\alpha(K(t_i,r_i,n_i))}{|K(t_i,r_i,n_i)|}: 1\leq i\leq k\right\}\prod_{1\leq
i\leq k}|K(t_i,r_i,n_i)|.$$
\end{theorem}

The circular graph $Circ(r,n)$ ($n\geq 2r$) is defined by: \[V(Circ(r,n))=\mathbb{Z}_n=\{0,1,\ldots,n-1\}\] and \[E(Circ(r,n))=\left\{(i,j):|i-j|\in \{r,r+1,\ldots,n-r\} \right\}.\]
It is well known that $\alpha(Circ(r,n))=r$. Mario and Juan \cite{MJ} determined the independence number of the direct products of circular graphs.
\begin{theorem}\label{cdp} (Mario and Juan \cite{MJ})
Let $n_i\geq 2r_i$ for $i=1,2\ldots,k$. Then
\[\alpha\left(\prod_{1\leq i\leq k}Circ(r_i,n_i)\right)=\max\left\{\frac{r_1}{n_1},\frac{r_2}{n_2},\ldots,\frac{r_k}{n_k}\right\}\prod_{1\leq i\leq k}n_i.\]
\end{theorem}

For positive integers $n$, let $S_n$ denote the permutation group on $[n]$. Two
permutations $f$ and $g$ are said to be intersecting if
there exists an $i\in [n]$ such that $f(i)=g(i)$.  We define a graph on
$S_{n}$ as that two  permutations are  adjacent if and
only if they are not intersecting. For brevity, this graph is also denoted by
$S_{n}$. Deza and
Frankl \cite{Deza} first obtained that $\alpha(S_n)= (n-1)!$.
 Cameron and Ku \cite{Cameron} proved that
each maximum
 independent set of $S_n$ is a coset of the stabilizer of a
 point, to which Larose and Malvenuto
\cite{Larose},  Wang and Zhang \cite{WJ} and Godsil and
 Meagher \cite{cgkm}  gave alternative proofs, respectively.
Recently, Cheng and Wong \cite{yku} further investigated the
independence number
 and the MIS-normality of the
 direct products of $S_n$.
 \begin{theorem}\label{si} (Cheng and Wong\cite{yku})
  Let $2\leq n_1=\cdots=n_p<n_{p+1}\leq\ldots,n_q$, $1\leq p\leq q$. Then
  $$\alpha\big(S_{n_1}\times S_{n_2}\times
  \cdots\times S_{n_q}\big)=(n_1-1)!
  \prod_{2\leq i\leq q}n_i!,$$
  and  the direct products $S_{n_1}\times S_{n_2}\times
  \cdots\times S_{n_q}$ is MIS-normal except for  the following cases:
\begin{enumerate}
  \item [\em(i)] $n_1=\cdots=n_p<n_{p+1}=3\leq n_{p+2}\leq \cdots\leq n_q$;
  \item [\em(ii)] $n_1=n_2=3\leq n_3\leq \cdots\leq n_q$;
  \item [\em(iii)] $n_1=n_2=n_3\leq n_4\leq \cdots\leq n_q$.
\end{enumerate}
\end{theorem}

In \cite{larose}, Larose and Tardif  investigated the relationship between projectivity and the structure of maximum
independent sets in powers of some vertex-transitive graphs, and obtained the MIS-normality of the powers of Kneser graphs and circular graphs.
\begin{theorem}(Larose and Tardif \cite{larose})
Let $n$ and $r$ be two positive integers. If $n>2r$, then  both $K^k(1,r,n)$ and $Circ^k(r,n)$ are MIS-normal for all positive integer $k$.
\end{theorem}

Besides the above results, Larose and Tardif \cite{larose} prove
that if  $G$ is vertex-transitive, then
$\alpha(G^n)=\alpha(G)|V(G)|^{n-1}$ for all $n>1$. They also ask
whether or  not $G^n$ is MIS-normal if $G^2$ is MIS-normal.
Recently, Ku and Mcmillan \cite{KM} gave an affirmative answer to
this problem, and we   solved this problem in a more general setting
\cite{zhhj}.

In  this paper we shall solve both Problem \ref{Tard} and Problem
\ref{MIS}. To state our results we need to introduce some notations
and notions.

For a graph $G$, let $I(G)$ denote the set of all maximum
independent sets of $G$. Given a subset  $A$  of $V(G)$, we define
$$N_G(A)=\{b\in V(G):\mbox{$(a,b)\in E(G)$ for some $a\in A$}\}
$$
$$N_G[A]=N_G(A)\cup A \mbox{ and } \bar{N}_G[A]=V(G)-N_G[A].$$
If $G$ is clear from the context, for simplicity, we will omit the
index $G$.

In \cite{zhhj}, by the so-called ``No-Homomorphism" lemma of
Albertson and Collins \cite{makl} we proved the following result.
\begin{proposition}\label{coro1}(\cite{zhhj})
Let  $G$ be a vertex-transitive graph. Then, for every independent
set $A$ of $G$, $\frac{|A|}{|N_G[A]|}\leq \frac{\alpha(G)}{|V(G)|}$.
Equality implies that  $|S\cap N_G[A]|=|A|$  for every $S\in
  I(G)$, and in particularly $A\subseteq S$ for some $S\in I(G)$.
\end{proposition}

An independent set $A$ in $G$ is said to be {\em imprimitive} if
$|A|<\alpha(G)$ and $\frac{|A|}{|N[A]|}=\frac{\alpha(G)}{|V(G)|}$.
And $G$ is called {\em IS-imprimitive} if $G$ has an imprimitive
independent set. In any other cases, $G$ is called
\emph{IS-primitive}. From definition we see that  a disconnected vertex-transitive graph $G$ is IS-imprimitive and hence an IS-primitive vertex-transitive graph $G$ is connected.

The following Theorem  is the main result of this paper.
\begin{theorem}\label{DP}
Let $G$ and $H$ be two  vertex-transitive graphs with
$\frac{\alpha(G)}{|G|}\geq\frac{\alpha(H)}{|H|}$. Then
\[\alpha(G\times H)=\alpha(G)|H|,\]
and either:
\begin{enumerate}
  \item [\rm(i)] $G\times H$ is MIS-normal, or
  \item [\rm(ii)] $\frac{\alpha(G)}{|G|}=\frac{\alpha(H)}{|H|}$ and one of them is IS-imprimitive, or
  \item [\rm(iii)]$\frac{\alpha(G)}{|G|}>\frac{\alpha(H)}{|H|}$ and $H$ is disconnected.
\end{enumerate}
\end{theorem}

  We leave the proof of  Theorem \ref{DP} to the next section, while in Section 3,
  we discuss the MIS-normality of the direct products of more than two
  vertex-transitive graphs.

\section{Proof of Theorem \ref{DP}}
  Let $S$ be a maximum independent set of $G\times
H$. Then  $|S|\geq \alpha(G)|H|\geq |G|\alpha(H)$. We now prove
$\alpha(G\times H)\leq \alpha(G)|H|$.

For every $a\in G$, define \[X_a=\{x\in H: (a,x)\in S\}.\] Since $S$
is an independent set of $G\times H$, for each $x\in X_a$ and $y\in
X_b$,  $(x,y)\not\in E(H)$ whenever  $(a,b)\in E(G)$. In this case,
we say that $X_a$ and $X_b$ are cross-independent. This concept is
equivalent to cross-intersecting families in extremal set theory. We
refer \cite{wz} for details.

In the language of cross-intersecting families, Borg
\cite{borg1,borg2,borg3} introduce a decomposition of $X_a$ as
follows.
\[
X^*_a=\{x\in X_a: N_H(x)\cap X_a=\emptyset\},\]\[ X'_a=\{x\in X_a:
N_H(x)\cap X_a\neq\emptyset\}\] and \[X'=\bigcup_{a\in V(G)}X'_a.\]
Clearly, $X_a^*$ is an independent set of $H$ for every  $a\in
V(G)$, and $|S|=\sum_{a\in V(G)}|X_a|$. Here, the empty set is
regarded as an independent set.

We list all distinct $X^*_a$'s  as $Y_1,Y_2,\ldots, Y_k$, and define
$$B_i=\{a\in V(G): X^*_a=Y_i\},\ i=1,2,\ldots,k.$$ We then obtain a
partition of $V(G)$ as $V(G)=B_{1}\cup B_{2}\cup\cdots\cup B_{k}$.
Then
\begin{eqnarray}\nonumber
|S|&=&\sum_{a\in V(G)}|X_a|=\sum_{a\in
V(G)}(|X_a^*|+|X'_a|)=\sum_{i=1}^k\ \sum_{a\in
B_{i}}|X_a^*|+\sum_{a\in V(G)}|X_a'|\\ &=&
\sum_{i=1}^k|Y_i||B_{i}|+\sum_{x\in X'}|A_x|,\label{IS1}
\end{eqnarray}
where
 \[A_x=\{a\in V(G): x\in X'_a\}.\]
For every pair $a,b\in V(G)$, it is easy to verify that
$(a,b)\not\in E(G)$  if $X'_a\cap X'_b\neq \emptyset$.
 Therefore, $A_x$ is an independent set of $G$.
By Proposition \ref{coro1} we have that
\begin{equation}\label{3}
|A_x|\leq \frac{\alpha(G)}{|V(G)|}|N_G[A_x]|,
\end{equation}
 and equality holds
if and only if $|A_x|=0$, or $|A_x|=\alpha(G)$, or $A_x$ is an
imprimitive independent set of $G$.

Suppose $x\in N_H[Y_i]=N_H(Y_i)\cup Y_i$.
 If
$x\in N_H(Y_i)$, then there exists $y\in Y_i$ such that $(x,y)\in
E(H)$ and $\{(a,x),(b,y)\}\subset S$ for any $b\in B_{i}$ and $a\in
A_x$, hence $(a,b)\not\in E(G)$ since $S$ is an independent
set; 
 if $x\in Y_i$, then for each $a\in A_x$, there is a $z\in X_a$ with
 $(x,z)\in E(H)$ and $\{(a,z),(b,x)\}\subset S$, yielding $(a,b)\not\in
E(G)$. Thus proving that $B_{i}\subseteq \bar{N}_G[A_x]$ if $x\in
N_H[Y_i]$. From this it follows that
$$\sum_{i: x\in N_H[Y_i]}|B_{i}|\leq
|\bar{N}_G[A_x]|=|V(G)|-|N_G[A_x]|,$$ i.e.,
\begin{equation}\label{NS}
|N_G[A_x]|\leq |V(G)|-\sum_{i: x\in N_H[Y_i]}|B_{i}|=\sum_{i: x\in
\bar N_H[Y_i]}|B_{i}|.
\end{equation}
Note that
\begin{equation}\label{X}
X'\subseteq\bigcup_{ i=1}^k\bar{N}_H[Y_i].
\end{equation}
Together with (\ref{3}), (\ref{NS})and (\ref{X}), we then obtain  that
\begin{eqnarray}
&&\sum_{x\in X'}|A_x|\leq \frac{\alpha(G)}{|V(G)|}\sum_{x\in
X'}\sum_{i:x\in\bar N_H[Y_i]}|B_{i}|\nonumber\\
&\leq& \frac{\alpha(G)}{|V(G)|}\sum_{i=1}^k\sum_{x\in\bar
N_H[Y_i]}|B_{i}|=\frac{\alpha(G)}{|V(G)|}\sum_{i=1}^k|B_{i}||\bar{N}_H[Y_i]|.\label{IS2}
\end{eqnarray}
Combining (\ref{IS1}) and (\ref{IS2}) gives that
\begin{eqnarray*}\label{8}
|S|&= &\sum_{i=1}^k|Y_i||B_{i}|+\sum_{x\in X'}|A_x|\\
&\leq&\sum_{i=1}^k|Y_i||B_{i}|+\frac{\alpha(G)}{|V(G)|}\sum_{i=1}^k|B_{i}||\bar{N}_H[Y_i]|\\
&=&\sum_{i=1}^k|B_{i}|\left(\frac{\alpha(G)}{|V(G)|}|H|+|Y_i|-\frac{\alpha(G)}{|V(G)|}|N_H[Y_i]|\right)\\
&=& \alpha(G)|H|+\sum_{i=1}^k|B_{i}|\left(|Y_i|-\frac{\alpha(G)}{|V(G)|}|N_H[Y_i]|\right)\\
 &\leq&
\alpha(G)|H|.
\end{eqnarray*}
The last inequality follows from  that
\begin{eqnarray}\label{le}|Y_i|-\frac{\alpha(G)}{|G|}|N_H[Y_i]|\leq
|Y_i|-\frac{\alpha(H)}{|V(H)|}|N_H[Y_i]|\leq 0,\end{eqnarray} by
Proposition \ref{coro1}.

The maximum of $|S|$ implies that $|S|=\alpha(G)|H|$, from which it
follows that  equalities (\ref{3}), (\ref{NS}), (\ref{X}) and (\ref{le}) hold.
Also, from Proposition \ref{coro1}, equality (\ref{le}) means that  either $Y_i=\emptyset$, or
$\frac{\alpha(G)}{|G|}=\frac{\alpha(H)}{|V(H)|}$ and  $Y_i$ is either imprimitive or a maximum independent set of $H$ for $i=1,2,\ldots, k$.

 We now prove that either $S$ is  the preimages
of projections of a maximum independent set of $G$ or $H$, or (ii) or (iii)
holds. There are two cases to be
considered.

Case 1: $\frac{\alpha(G)}{|G|}>
\frac{\alpha(H)}{|H|}$. Then, equality (\ref{le}) means that $Y_i=\emptyset$ for all $i$, and so $X'=V(H)$ by equality (\ref{X}). Hence, from equality (\ref{3}) it follows
 that $A_x$ is a maximum independent set of $G$ for all $x\in V(H)$.  With this assumption we have that for any
$x,y\in V(H)$ with $(x,y)\in E(H)$, if $A_x\neq A_y$, there must
exist $a\in A_x$ and $b\in A_y$ with $(a,b)\in E(G)$ since both
$A_x$ and $A_y$ are maximum independent set, so $[(a,x),(b,y)]\in
E(G\times H)$, contradicting
 $\{(a,x),(b,y)\}\subset S$. Therefore, $A_x=A_y$ whenever
 $(x,y)\in E(H)$, which implies that $S$ is the
preimage of a maximum independent set of $G$ under projections  if $H$ is
connected.

Case 2: $\frac{\alpha(G)}{|G|}=
\frac{\alpha(H)}{|H|}$. Then, equality (\ref{le}) means that either $|Y_i|=0$ or $\alpha(H)$,
or $Y_i$ is an imprimitive independent set of $H$ for each index $i$. If $Y_i$ is an imprimitive independent
set of $H$ for some $i$, then $H$ is IS-imprimitive. If $|Y_i|=\alpha(H)$ for all $i$, then $X_a=X_a^*$
is a maximum independent set of $H$ for all $a\in V(G)$, and we can  prove in the similar way as in Case 1
that $S$ is the preimage of a maximum independent set of $H$ under projections if $G$ is connected.
 We now suppose that $|Y_i|=0$  for some $i$. With this assumption, then equality (\ref{X}) implies $X'=V(H)$,
  and then equality (\ref{NS}) means that either $A_x$ is either imprimitive or a maximum independent set of $G$
  for all $x\in V(H)$. If the former holds for some $x\in V(H)$, we  have that $H$ is IS-imprimitive; otherwise,
  the latter holds for all $x\in V(H)$, and then we can  prove in the similar way as in Case 1 that $S$ is the
  preimage of a maximum independent set of $G$ under projections if $H$ is connected.

\section{Concluding Remark.}

Let $G_1,G_2, \ldots,G_n$ be $n$ non-empty vertex-transitive graphs,
and  set $G=G_1\times G_2\times \cdots\times G_n$. From Theorem
\ref{DP} it immediately follows that\[
\alpha(G)=\alpha(G_1)\prod_{2\leq i\leq n}|G_i|.\] We now discuss
the MIS-normality of $G$. For convenience, we say $G$ is MIS-normal
if $n=1$.

 A graph $H$ is said to be \textit{non-empty} if $E(H)\neq
\emptyset$. It is well known that if $H$ is a non-empty
vertex-transitive graph,  then $\frac{\alpha(H)}{|H|}\leq \frac 12$,
and equality holds if and only if $H$ is a bipartite graph.

Without loss of generality we may assume that $\frac 12\geq
\frac{\alpha(G_1)}{|G_1|}=\cdots=\frac{\alpha(G_\ell)}{|G_\ell|}>\frac{\alpha(G_{\ell+1})}{|G_{\ell+1}|}
\geq \cdots \geq \frac{\alpha(G_n)}{|G_n|}$, and write
$H_0=G_1\times \cdots\times G_\ell$ and $H_i=H_{i-1}\times
G_{\ell+i}$ for $i=1,\ldots,n-\ell$ subject to $n>\ell$. Then
$G=H_{n-\ell}$ and with
$\frac{\alpha(H_{i-1})}{|H_{i-1}|}>\frac{\alpha(G_{\ell+i})}{|G_{\ell+i}|}$ for
$i\geq 1$.

\begin{proposition}\label{p1} Suppose $n>\ell$. Then $G$ is MIS-normal if and
only if $H_0$ is MIS-normal and $G_{\ell+1},\ldots, G_n$ are all
connected.
\end{proposition}
\noindent{\bf Proof}. Since
$\alpha(G)=\alpha(H_0)\prod_{i=\ell+1}^{n}|G_i|$, we have that if
$H_0$ is not MIS-normal, then  $G$ is not MIS-normal. Furthermore, if
$G_i$ is not connected for for some $i\geq 1$, writing $G_i=G_i'\cup
G_i''$, a union of disjoint subgraphs,
  then,  for all $I_1,I_2\in I(H_{i-1})$ with $I_1\neq I_2$, it is clear that
  $S=(I_1\times G_i')\cup (I_2\times G_i'')\in I(H_i)$, which
  is not a preimage of any independent set of one factor under
  projections, i.e., $H_i$ is not MIS-normal, hence $G$ is not MIS-normal.

 Conversely,  suppose  $H_{0}$ is MIS-normal, and $G_{\ell+i}$ is connected for $i\geq
 1$. Since
 $\frac{\alpha(H_{i-1})}{|H_{i-1}|}>\frac{\alpha(G_{\ell+i})}{|G_{\ell+i}|}$,  Theorem \ref{DP} implies that
each maximal-sized independent set is of the form $S\times
G_{\ell+i}$, where $S\in I(H_{i-1})$, which means that $H_i$ is
MIS-normal for $i\geq 1$. We thus prove that $G$ is MIS-normal. \qed

We now discuss the case $n=\ell$, that is, each $G_i$ has the
identical independence ratio. To deal with this case we need a lemma
as follows.
\begin{lemma}\label{bipart} Suppose that G is a vertex-transitive bipartite graph. Then $G$ is
imprimitive if and only if $G$ is disconnected.
\end{lemma}
\textbf{Proof.} It is clear that $G$ is imprimitive if $G$ is
disconnected. On the converse, if $G$ is imprimitive, then there is
an imprimitive independent set $A$ such that
$\frac{|A|}{|N_G[A]|}=\frac{\alpha(G)}{|G|}=\frac{1}{2}.$ Set
$B=N_G(A)$.  $|B|=|A|$ and $A\subseteq N_G(B)$ is clearly. If
$N_G(B)\neq A$, then we obtain that $\sum_{u\in A}d(u)\leq
\sum_{v\in B}d(v)$, which induces a contradiction. Hence $N_G(B)=A$,
that is to say $G$ is disconnected.
\begin{proposition}\label{p2} Suppose that $
\frac{\alpha(G_1)}{|G_1|}=\cdots=\frac{\alpha(G_n)}{|G_n|}=\frac{\alpha(G)}{|G|}$.
Then $G$ is MIS-normal if and only if one of the following holds.

{\rm (i)}  $\frac{\alpha(G)}{|G|}<\frac 12$ and every $G_i$ is
IS-primitive.

{\rm (ii)} $\frac{\alpha(G_1)}{|G_1|}=\frac 12$,   $n=2$ and both
$G_1$ and $G_2$ are connected.
\end{proposition}
\noindent{\bf Proof}. For $1\leq i\leq n$, set $\hat
G_i=G_1\times\cdots\times G_{i-1}\times G_{i+1}\times\cdots\times
G_n$. Then $G=\hat G_i\times G_i$ for $i=1,2,\ldots,n$. If $G_i$ is
imprimitive, letting $A_i$ be an imprimitive independent set of
$G_i$, for every $I\in I(\hat G_i)$, it is easy to see that
 $S=(\hat G_i\times A_i)\cup (I\times \bar N_{G_i}[A_i])\in I(G)$, which
  is not a preimage of any independent set of $\hat G_i$ or $G_i$ under
  projections, therefore, $G$ is not MIS-normal. Conversely, if both
  $\hat G_i$ and $G_i$ are IS-primitive, Theorem \ref{DP} implies
  that $G$ is MIS-normal. It remains to check when $\hat G_i$ is
  IS-primitive. Summing up the above, $G$ is MIS-normal if and only
  if both $\hat G_i$ and $G_i$ are IS-primitive. To complete the
  proof, it remains to check when $\hat G_i$  is IS-primitive.
We distinguish  two cases.

  Case (i):$\frac{\alpha(G)}{|G|}<\frac 12$. In this case,
  Theorem 2.6 in \cite{zhhj} says that if $G$ is MIS-normal, then
  both $\hat G_i$ and $G_i$ are IS-primitive. The induction implies (i).

  Case (ii): $\frac{\alpha(G_1)}{|G_1|}=\frac 12$, i.e., every  $G_i$ is
  bipartite. From Lemma \ref{bipart} it follows that $\hat G_i$ and $G_i$ is IS-primitive  if and only if
  both $\hat G_i$ and $G_i$ are connected. However, it is well known
  that $\hat G_i$ is disconnected if $n>2$, thus proving (ii).
  \qed

Combining Proposition \ref{p1} and Proposition \ref{p2} gives the
following theorem.

\begin{theorem} Let $G_1,G_2, \ldots,G_n$ be  connected vertex-transitive graphs
with $\frac 12\geq
\frac{\alpha(G_1)}{|G_1|}=\cdots=\frac{\alpha(G_\ell)}{|G_\ell|}>\frac{\alpha(G_{\ell+1})}{|G_{\ell+1}|}
\geq \cdots \geq \frac{\alpha(G_n)}{|G_n|}$, where  $n\geq 2$ and
$1\leq \ell\leq n$. Then $G_1\times G_2\times\cdots\times G_n$ is
MIS-normal if and only if one of the following holds:

{\rm (i)} $\frac{\alpha(G_1)}{|G_1|}<\frac 12$ and $G_1,G_2,
\ldots,G_\ell$ are all IS-primitive whenever $\ell>1$.

{\rm (ii)} $\frac{\alpha(G_1)}{|G_1|}=\frac 12$ and $\ell\leq 2$.
\end{theorem}

\noindent\textbf{Acknowledgement}
  The author is greatly indebted to Professor J. Wang for giving
  useful comments, suggestions and helps that have considerably improved the
  manuscript.

\end{document}